\documentclass[12pt, reqno]{amsart}
\usepackage{amsmath, amsthm, amscd, amsfonts, amssymb, graphicx, color}
\usepackage[bookmarksnumbered, colorlinks, plainpages]{hyperref}

\newtheorem{theorem}{Theorem}[section]

\theoremstyle{definition}
\newtheorem{definition}[theorem]{Definition}

\theoremstyle{remark}

\numberwithin{equation}{section}

\begin{document}
\setcounter{page}{1}


\vspace{1.3 cm}
\title[n-multipliers]{n-multipliers and their relations with n-homomorphisms}

\author[J. Laali, M. Fozouni]{Javad Laali and Mohammad fozouni}

\address{Department of Mathematics, Faculty of mathematics and computer sciences, Kharazmi university, 43 Mofateh Ave, Tehran, Iran.}
\email{\textcolor[rgb]{0.00,0.00,0.84}{Lali@khu.ac.ir}}
\address{Department of Mathematics, Faculty of Sciences and Engineering, Gonbad Kavous University, P. O. Box 163, Gonbad-e Kavous, Golestan, Iran}
\email{\textcolor[rgb]{0.00,0.00,0.84}{fozouni@gonbad.ac.ir}}


\subjclass[2010]{Primary: 46H05, Secondary: 42A45.}

\keywords{ Banach algebra, multiplier, tensor product space, Banach module.}


\begin{abstract}

Let $A$ be a Banach algebra and $X$ be a Banach $A$-bimodule. We introduce and study the notions of $n$-multipliers and approximate local $n$-multipliers by generalizing the classical concept of multipliers from $A$ into $X$. As an algebraic result, we construct a  Banach algebra consisting of $n$-multipliers on $A$ and under some mild conditions, we give a nice relation of this algebra with  $n$-homomorphisms from $A$ into $\mathbb{C}$.
\end{abstract}
\maketitle
\section{Introduction and preliminaries}
\noindent The concept of a multiplier first appears in harmonic analysis in connection with the
theory of summability for Fourier series. Subsequently, the notion has been employed in other areas of harmonic analysis, such as the investigation of homomorphisms of group algebras, in the general theory of Banach algebras, and so on; see \cite{Larsen}. Many authors generalized the notion of a multiplier in different ways. See \cite{McKennon} and \cite{Adib}, for one of this generalizations.

In this paper, our main concern will not be with these applications of the theory of multipliers and its generalizations. We only develop the theory of multipliers  differently from the previous ways,
by introducing a new class of operators from a Banach algebra $A$ into a Banach $A$-bimodule $X$.

Let $A$ be a Banach algebra and $a, b\in A$. Define a bounded bilinear functional on $A^{*}\times A^{*}$ as
\begin{center}
$(a\otimes b)(f,g)=f(a)g(b)\hspace{1cm}(f, g\in A^{*})\cdot$
\end{center}

The projective tensor product space $A\widehat{\otimes}A$ is a Banach algebra and a Banach $A$-bimodule that is characterized as follows
$$
\left\{\Sigma_{n=1}^{\infty}a_{n}\otimes b_{n}: n\in \mathbb{N}, a_{n},b_{n}\in A, \Sigma_{n=1}^{\infty}||a_{n}||||b_{n}||<\infty\right\},
$$
and its module actions are defined by
\begin{center}
$a.(b\otimes c)=ab\otimes c,\hspace{0.5cm} (b\otimes c).a=b\otimes ca\hspace{1cm}(a,b,c\in A)\cdot$
\end{center}
A Banach algebra $A$ is called \emph{nilpotent} if there exists an integer  $n\geq 2$  such that
\begin{center}
$A^{n}=\{a_{1}a_{2}a_{3}...a_{n}: a_{1}, a_{2}, a_{3},...,a_{n}\in A\}=\{0\}\cdot$
\end{center}
The minimum of numbers $n$ that $A^{n}=\{0\}$ is called the $index$ of $A$ which we denote it by $I(A)$, i.e., if $I(A)=n$, then there exists $a_1, a_2,\ldots, a_{n-1}\in A$ such that $a_1a_2\ldots a_{n-1}\neq 0$.

 To see an  example of a nilpotent Banach algebra, suppose that $B$ is a Banach algebra and let $A$ be defined as follows
\begin{equation*}
A=\begin{bmatrix}
0 & B & B & B & B\\
0 &  0 & B & B & B\\
0 & 0 &  0 & B & B\\
0 & 0 & 0 & 0 & B\\
0 & 0 & 0 & 0 & 0\\
\end{bmatrix}.
\end{equation*}
Then, $A$ is a Banach algebra equipped with the usual matrix-like operations and $l_{\infty}$-norm such that $A$ is nilpotent with $I(A)=5$.

For undefined concepts and notations appearing in the sequel, one can consult \cite{Dales}.
\section{n-multipliers}
We start this section with the main object of the paper.
\begin{definition} Let $A$ be a Banach algebra, $X$ be a Banach $A$-bimodule and $T:A\rightarrow X$ be a bounded linear map. We say that $T$ is an \emph{n-multiplier}$ (n\geq 2)$ if
\begin{equation*}
T(a_{1}a_{2}...a_{n})=a_{1}\cdot T(a_{2}...a_{n})\quad(a_{1}, a_{2}, a_{3},...,a_{n}\in A)\cdot
\end{equation*}
\end{definition}
We will denote by $\mathrm{Mul}_{n}(A,X)$ the set of all $n$-multipliers of Banach algebra $A$ into $X$.
Now, we study in more details the space $\mathrm{Mul}_{n}(A,X)$ when $n\geq 3$ (in the case $n=2$ this is the space of all multipliers in the classical sense).

Let $A$ be a Banach algebra and $X$ be a Banach $A$-bimodule. The set $\mathrm{Mul}_{n}(A,X)$ is a vector subspace of $B(A,X)$; the space of all bounded linear maps from $A$ into $X$.

As the first result we show that $\mathrm{Mul}_{n}(A,X)$ is a closed vector subspace of $B(A,X)$.

\begin{theorem}\label{th2} Let $A$ be a Banach algebra and $X$ be a Banach $A$-bimodule. Then for all integers $n\geq 3$, the space $\mathrm{Mul}_{n}(A,X)$ is a closed vector subspace of $B(A,X)$.
\end{theorem}
\begin{proof}
We claim that $\mathrm{Mul}_{n}(A,X)$ is closed in $B(A,X)$. Suppose that $\{T_{m}\}$ is a sequence in $\mathrm{Mul}_{n}(A,X)$  such that converges to $T\in B(A,X)$.

Let $a_{1}, a_{2},..., a_{n}$ be arbitrary elements of $A$. So, we have
\begin{align*}
||T(a_{1}...a_{n})-a_{1}\cdot T(a_{2}...a_{n})||&\leq ||T(a_{1}...a_{n})-T_{m}(a_{1}...a_{n})||\\
&+||T_{m}(a_{1}...a_{n})-a_{1}\cdot T(a_{2}...a_{n})||\\
&\leq ||T-T_{m}||||a_{1}...a_{n}||\\
&+ ||a_{1}\cdot T_{m}(a_{2}...a_{n})-a_{1}\cdot T(a_{2}...a_{n})||\\
&\leq ||T-T_{m}||||a_{1}...a_{n}||+||T-T_{m}||||a_{1}||||a_{2}...a_{n}||\cdot
\end{align*}
If $m\rightarrow\infty$, we conclude that $T(a_{1}...a_{n})=a_{1}\cdot T(a_{2}...a_{n})$. The rest of the proof is easy.
\end{proof}

Similarly, one can see that $\mathrm{Mul}_{n}(A,X)$ is complete in the \emph{strong operator topology }(SOT), i.e., in the topology on $B(A,X)$ for which a net $\{T_{\alpha}\}$ converges to $T$  if and only if for each $a\in A$, $||Ta-T_{\alpha}a||\rightarrow 0$.

In the next two theorems, we give some relations between  the spaces of $n$-multipliers.

\begin{theorem} There exists a Banach algebra A and a Banach A-bimodule X such that for all positive integers $n\geq 3$
\begin{equation*}
\mathrm{Mul}_{2}(A,X)\subsetneq \mathrm{Mul}_{3}(A,X)\subsetneq\ldots \subsetneq \mathrm{Mul}_{n}(A,X)\cdot
\end{equation*}
\end{theorem}
\begin{proof} For every positive integer $n\geq 3$ take $A$ a nilpotent Banach algebra with $I(A)=n$ and $X=A\widehat{\otimes}A$. So, there exists non-zero elements $a_{1}, a_{2},...a_{n-1}\in A$ such that $a_{1}a_{2}...a_{n-1}\neq 0$.

The verification of the above chain of inclusions is easy. We only show each of the strict relations. For every integer number $i$ such that $2\leq i < n$, define a linear map $T_{i}:A\rightarrow A\widehat{\otimes}A$  by
\begin{equation*}
T_{i}(a)=a_{1}a_{2}...a_{n-1}\otimes a_{1}a_{2}...a_{n-(i+1)}a\hspace{1cm}(a\in A)\cdot
\end{equation*}
Thus, $T_{i}$ is an element of $\mathrm{Mul}_{i+1}(A, A\widehat{\otimes}A)$, but it does not belong to $\mathrm{Mul}_{i}(A, A\widehat{\otimes}A)$. To see this, let $f\in A^{*}$ be a functional such that $f(a_{1}a_{2}\ldots a_{n-1})\neq 0$. So
\begin{align*}
T_{i}(a_{n-i}\ldots a_{n-1})(f,f)&=(a_{1}a_{2}\ldots a_{n-1}\otimes a_{1}a_{2}\ldots a_{n-1})(f,f)\\
&=f(a_{1}a_{2}\ldots a_{n-1})^{2}\neq 0\cdot
\end{align*}
Therefore, $T_{i}(a_{n-i}\ldots a_{n-1})\neq 0$, but $a_{n-i}.T_{i}(a_{(n-i)+1}...a_{n-1})=0$ and this completes the proof.
\end{proof}


For a Banach algebra $A$, let $A^{2}=\textmd{span}\{ab:a,b \in A\}$. The Banach algebra $A$ is \emph{essential} if $\overline{A^{2}}=A$.

\begin{theorem}\label{Th: essential} Let $A$ be an essential Banach algebra  and $X$ be a Banach $A$-bimodule. Then for all integers $n\geq 3$ we have
\begin{equation*}
\mathrm{Mul}_{n-1}(A,X)=\mathrm{Mul}_{n}(A,X)\cdot
\end{equation*}
Specially, $\mathrm{Mul}_{2}(A,X)=\mathrm{Mul}_{n}(A,X)$ for all $n\geq 2$.
\end{theorem}

\begin{proof} Let $T\in \mathrm{Mul}_{n}(A,X)$, $A$ be essential and $a_{1}, ..., a_{n-1}$ be arbitrary elements of $A$. We show that $T\in \mathrm{Mul}_{n-1}(A,X)$. Since $a_{2}\in A=\overline{A^{2}}$, there exists a net $\{a_{2,\alpha}\}\in A^{2}$ with $a_{2,\alpha}=\sum_{i_{\alpha}} \beta_{i_{\alpha}}b_{i_{\alpha}}c_{i_{\alpha}}$, for, $b_{i_{\alpha}}, c_{i_{\alpha}}\in A$ and $\beta_{i_{\alpha}}\in \mathbb{C}$, such that $a_{2}=\lim_{\alpha}a_{2,\alpha}$.  So, we have

\begin{align*}
T(a_{1}a_{2}...a_{n-1})&=\lim_{\alpha}T(a_{1}(\sum_{i_{\alpha}} \beta_{i_{\alpha}}b_{i_{\alpha}}c_{i_{\alpha}})a_{3}...a_{n-1})\\
&=\lim_{\alpha}\sum_{i_{\alpha}}\beta_{i_{\alpha}}T(\overbrace{a_{1}b_{i_{\alpha}}c_{i_{\alpha}}a_{3}...a_{n-1}}^{n})\\
&=\lim_{\alpha}\sum_{i_{\alpha}}\beta_{i_{\alpha}}a_{1}\cdot T(b_{i_{\alpha}}c_{i_{\alpha}}a_{3}...a_{n-1})\\
&=a_{1}\cdot T(a_{2}...a_{n-1}),
\end{align*}
which completes the proof.
\end{proof}
\section{Relations with n-homomorphisms}

Let $A$ be a Banach algebra and $n\geq 3$ be a positive integer. Here we show the space $\mathrm{\mathrm{Mul}}_{n}(A,A)$ briefly by $\mathrm{Mul}_{n}(A)$.

We know that the space of all multipliers on $A$ is a Banach subalgebra of $B(A)=B(A,A)$ with composition of operators as product and the operator norm. But in general the space of $n$-multipliers on $A$ is not an algebra with composition of operators. So, we should define another product on this space to make $\mathrm{Mul}_{n}(A)$ into a Banach algebra.

Now, let $a_{0}\in A$ and consider $\bullet_{a_{0}}:\mathrm{Mul}_{n}(A)\times \mathrm{Mul}_{n}(A)\rightarrow \mathrm{Mul}_{n}(A)$ which is defined by
\begin{equation}\label{product def}
S\bullet_{a_{0}} T(a):=S(T(a)a_{0}^{n-2})\hspace{0.5cm}(a\in A),
\end{equation}
Without losing the generality we assume that $||a_{0}||\leq 1$.

\begin{theorem} Let $A$ be a Banach algebra. Then for all positive integers $n\geq 3$, $\mathrm{Mul}_{n}(A)$ is a Banach algebra, with the product $"\bullet_{a_{0}}"$ and  the operator norm.
\end{theorem}
\begin{proof} Clearly, $\mathrm{Mul}_{n}(A)$ is a vector space  with operations that inherit from $B(A)$. Let $S, T\in \mathrm{Mul}_{n}(A)$. First we show that $S\bullet_{a_{0}}T$ is well-defined, i.e., $S\bullet_{a_{0}}T\in \mathrm{Mul}_{n}(A)$. Let $a_{1}, a_{2},...,a_{n}\in A$, we have
\begin{align*}
S\bullet_{a_{0}}T(a_{1}...a_{n})=S(T(a_{1}...a_{n})a_{0}^{n-2})&=S(a_{1}T(a_{2}...a_{n})a_{0}^{n-2})\\
&=a_{1} S(T(a_{2}...a_{n})a_{0}^{n-2})\\
&=a_{1}(S\bullet_{a_{0}}T)(a_{2}...a_{n})\cdot
\end{align*}
Therefore,  $S\bullet_{a_{0}}T$ is an $n$-multiplier on $A$.

Let $T_{1}, T_{2}$ and $T_{3}$ be elements of $\mathrm{Mul}_{n}(A)$. We have
\begin{align*}
&(T_{1}\bullet_{a_{0}} T_{2})\bullet_{a_{0}} T_{3}(a)=(T_{1}\bullet_{a_{0}} T_{2})(T_{3}(a)a_{0}^{n-2})=T_{1}(T_{2}(T_{3}(a)a_{0}^{n-2})a_{0}^{n-2})\\
&T_{1}\bullet_{a_{0}} (T_{2}\bullet_{a_{0}} T_{3})(a)=T_{1}((T_{2}\bullet_{a_{0}}T_{3})(a)a_{0}^{n-2})=T_{1}(T_{2}(T_{3}(a)a_{0}^{n-2})a_{0}^{n-2})
\end{align*}
Hence, the product $"\bullet_{a_{0}}"$ is  associative.

On the other hand, Theorem \ref{th2} shows that $\mathrm{Mul}_{n}(A)$ is a closed vector subspace of $B(A)$ with the operator norm. The investigation of the other properties are easy.
\end{proof}
Recall that for Banach algebras $A$ and $B$ a linear map $\phi:A\rightarrow B$ is called an \emph{n-homomorphism} if, $\phi(a_{1}a_{2}\ldots a_{n})=\phi(a_{1})\phi(a_{2})\ldots\phi(a_{n})$  for all $a_{1}, a_{2}, \ldots, a_{n}\in A$ \cite{Hejazian}.

For each integer $n\geq 2$, suppose that $\Delta_{n}(A)$ denotes the \emph{n-character space} of $A$, i.e., the space consisting of all non-zero $n$-homomorphisms from $A$ into $\mathbb{C}$. It is clear that for every integer $n\geq 3$, $\Delta_{2}(A)\subseteq \Delta_{n}(A)$. The last inclusion may be strict. As an example for $n=3$, if $\phi\in \Delta_{2}(A)$, then $\varphi:=-\phi$ is in $\Delta_{3}(A)$, but $\varphi$ is not a 2-character from $A$ into $\mathbb{C}$.

Let $\phi_{n}\in \Delta_{n}(A)$. Define $\widetilde{\phi}_{n}:(\mathrm{Mul}_{n}(A), \bullet_{a_{0}}) \rightarrow \mathbb{C}$ by
\begin{equation*}
\widetilde{\phi}_{n}(T)=\phi_{n}(T(a_{0}^{n-1}))\quad (T\in \mathrm{Mul}_{n}(A))\cdot
\end{equation*}
Clearly, $\widetilde{\phi}_{n}$ is a linear operator. We say that $\widetilde{\phi}_{n}$ extends $\phi_{n}$ if, $\widetilde{\phi}_{n}(L_{a})=\phi_{n}(a)$ for all $a\in A$.

In the next theorem, under some mild conditions, we show that $\Delta_{n}(\mathrm{Mul}_{n}(A))\neq \emptyset$. Recall that $Z(A)$ denotes the center of $A$, i.e.,
 \begin{center}
 $Z(A)=\{a\in A : ab=ba\ (b\in A)\}$.
 \end{center}
\begin{theorem}\label{Th: n-hom} Let $A$ be a Banach algebra and let $a_0\in Z(A)\setminus\{0\}$. Then $\widetilde{\phi}_{n}\in \Delta_{n}(\mathrm{Mul}_{n}(A))$ is an extension of $\phi_{n}\in \Delta_{n}(A)$ if $\phi_{n}(a_{0})=1$.
\end{theorem}
\begin{proof} Suppose that there exists $\phi_{n}\in \Delta_{n}(A)$ with $\phi_{n}(a_{0})=1$. We must show that $\widetilde{\phi}_{n}$ is a non-zero $n$-homomorphism. For each $a\in A$ we have
\begin{equation*}
\widetilde{\phi}_{n}(L_{a})=\phi_{n}(L_{a}(a_{0}^{n-1}))=\phi_{n}(aa_{0}^{n-1})=\phi_{n}(a)\phi_{n}(a_{0})^{n-1}=\phi_{n}(a)\cdot
\end{equation*}
Therefore, in the especial case when $a=a_{0}$, we have $\widetilde{\phi}_{n}(L_{a_{0}})=1$. So, $\widetilde{\phi}_{n}$ is a non-zero extension of $\phi_{n}$.

On the other hand, for $T_{1}, T_{2}, \ldots, T_{n}\in \mathrm{Mul}_{n}(A)$, we have
\begin{align*}
\widetilde{\phi}_{n}(T_{1}\bullet_{a_{0}}\ldots\bullet_{a_{0}}T_{n})&=\phi_{n}(T_{1}\bullet_{a_{0}}\ldots\bullet_{a_{0}}T_{n}(a_{0}^{n-1}))\\
&=\phi_{n}(T_{1}(T_{2}\overbrace{(\ldots(}^{n-2} T_{n}(a_{0}^{n-1})a_{0}^{n-2}\overbrace{)\ldots)}^{n-2}a_{0}^{n-2}))\\
&=\phi_n(a_{0})^{n-1}\phi_{n}(T_{1}(T_{2}\overbrace{(\ldots(}^{n-2} T_{n}(a_{0}^{n-1})a_{0}^{n-2}\overbrace{)\ldots)}^{n-2}a_{0}^{n-2}))\\
&=\phi_{n}(a_{0}^{n-1}T_{1}(T_{2}\overbrace{(\ldots(}^{n-2} T_{n}(a_{0}^{n-1})a_{0}^{n-2}\overbrace{)\ldots)}^{n-2}a_{0}^{n-2}))\\
&=\phi_{n}(T_{1}(T_{2}\overbrace{(\ldots(}^{n-2} T_{n}(a_{0}^{n-1})a_{0}^{n-1}\overbrace{)\ldots)}^{n-2}a_{0}^{n-1}))\\
&\ \ \ \ \ \ \ \ \ \ \ \ \ \ \ \ \ \ \quad\quad\quad\quad\vdots\\
&=\phi_{n}(T_{1}(a_{0}^{n-1})T_{2}(a_{0}^{n-1})\ldots T_{n}(a_{0}^{n-1}))\\
&=\phi_{n}(T_{1}(a_{0}^{n-1}))\phi_{n}(T_{2}(a_{0}^{n-1}))\ldots \phi_{n}(T_{n}(a_{0}^{n-1}))\\
&=\widetilde{\phi}_{n}(T_{1})\widetilde{\phi}_{n}(T_{2})\ldots \widetilde{\phi}_{n}(T_{n})\cdot
\end{align*}
Therefore, $\widetilde{\phi}_{n}$ is an $n$-homomorphism and this completes the proof.
\end{proof}
\section{Approximate local  n-multipliers}

In \cite{Samei}, Samei investigated the approximate local left 2-multipliers and study some of its relations with left $2$-multipliers on a Banach algebra $A$. In this section we give two theorems similar to  Theorem 2.2 and Proposition 2.3 of \cite{Samei}  for  $n$-multipliers. Indeed,  we are interested in determining when an approximately local  $n$-multiplier (Definition \ref{Def}) is an $n$-multiplier. First we give the following definition.

\begin{definition}\label{Def} Let $X$ be a Banach $A$-module and $T:A\rightarrow X$ be a bounded linear operator. We say that $T$ is an \emph{approximately local  n-multiplier} if, for each $a\in A$ there exists a sequence $\{T_{a,m}\}$ of $n$-multipliers such that $T(a)=\lim_{m}T_{a,m}(a)$.
\end{definition}

We recall the algebraic reflexivity from  \cite{Conway}. 
Let $X$ and $Y$ be Banach spaces and $S$ be a subset of $B(X,Y)$. Put
\begin{equation*}
\textmd{ref}(S)=\{T\in B(X,Y): T(x)\in \overline{\{s(x):s\in S\}}\ \ (x\in X)\}.
\end{equation*}
Then  $S$ is algebraically reflexive if, $S=\textmd{ref}(S)$ or just $\textmd{ref}(S)\subseteq S$.

\begin{theorem}\label{th4} Let $A$ be a Banach algebra and $X$ be a  Banach $A$-module. Then the following statements are equivalent.
\begin{enumerate}
  \item Every approximately local n-multiplier from $A$ into $X$ is an  n-multiplier $(n\geq 3)$.
  \item $\mathrm{Mul}_{n}$(A,X) is algebraically reflexive.
\end{enumerate}
\end{theorem}
\begin{proof} $(1)\rightarrow (2)$: Let $T\in \textmd{ref}(\mathrm{Mul}_{n}(A,X))$. So, for all $a\in A$ there exists a sequence $\{T_{m}\}$ in $\mathrm{Mul}_{n}(A,X)$ such that, $T(a)=\lim_{m}T_{a,m}(a)$. Hence, $T$ is an  approximately local  $n$-multiplier. Therefore,  $T$ is an  $n$-multiplier by assumption and this shows that $\mathrm{Mul}_{n}$(A,X) is algebraically reflexive.

$(2)\rightarrow (1)$: Let $T:A\rightarrow X$ be an approximately local  $n$-multiplier. So, for all $a\in A$, there exists a sequence $\{T_{a,m}\}$ such that, $T(a)=\lim_{m}T_{a,m}(a)$. Hence, $T\in \textmd{ref}(\mathrm{Mul}_{n}(A,X))$ and reflexivity of $\mathrm{Mul}_{n}(A,X)$ implies that $T$ is an $n$-multiplier.
\end{proof}

Let $A$ be a Banach algebra and $X$ be a  Banach $A$-module. Then for each $x\in X$, the \emph{left annihilator} of $x$ in $A$ is defined by $x^{\bot}=\{a\in A: a\cdot x=0\}$.

\begin{theorem}\label{th5} Suppose that $A$ is a Banach algebra such that $\mathrm{Mul}_{n}(A,A^{*})$ is algebraically reflexive and $X$ is a left Banach $A$-module with $\{x\in X: x^{\bot}=A\}=0$. Then every approximately local  n-multiplier from $A$ into $X$ is an n-multiplier.
\end{theorem}
\begin{proof} Let $T:A\rightarrow X$ be an approximate local $n$-multiplier and $f\in X^{*}$. Define a map  $\mathfrak{M}_{f}:X\rightarrow A^{*}$ as follows
\begin{equation*}
\mathfrak{M}_{f}(x)=x\bullet f\hspace{0.5cm}(x\in X)
\end{equation*}
where $x\bullet f\in A^{*}$ is defined by $x\bullet f(a)=f(a\cdot x)$ for all $a\in A$.
Therefore, $\mathfrak{M}_{f}$ is a bounded right $A$-module morphism. Because, for $a\in A$ and $x\in X$ we have
\begin{equation*}
\mathfrak{M}_{f}(a\cdot x)=(a\cdot x)\bullet f=a\cdot (x\bullet f)=a\cdot \mathfrak{M}_{f}(x).
\end{equation*}
So
\begin{equation*}
\mathfrak{M}_{f}\circ T\in \textmd{ref}(\mathrm{Mul}_{n}(A,A^{*}))=\mathrm{Mul}_{n}(A,A^{*}).
\end{equation*}
Thus, $\mathfrak{M}_{f}\circ T\in \mathrm{Mul}_{n}(A,A^{*})$. Now, for  $a_{1}, a_{2}, a_{3},... a_{n}\in A$ we have
\begin{align*}
\mathfrak{M}_{f}(T(a_{1}a_{2}a_{3}... a_{n}))=\mathfrak{M}_{f}\circ T(a_{1}a_{2}a_{3}... a_{n})
&=a_{1}\cdot\mathfrak{M}_{f}\circ T(a_{2}a_{3}... a_{n})\\
&=a_{1}\cdot\mathfrak{M}_{f}(T(a_{2}a_{3}... a_{n}))\\
&=\mathfrak{M}_{f}(a_{1}\cdot T(a_{2}a_{3}... a_{n})).
\end{align*}
Therefore, $\mathfrak{M}_{f}(T(a_{1}a_{2}a_{3}... a_{n})-a_{1}\cdot T(a_{2}a_{3}... a_{n}))=0$. If we put
\begin{equation*}
u=T(a_{1}a_{2}a_{3}... a_{n})-a_{1}\cdot T(a_{2}a_{3}... a_{n}),
\end{equation*}
then $f(a\cdot u)=0$ for all $a\in A$. So, by Hahn-Banach's theorem we have $a\cdot u=0$ for all $a\in A$. So, $u^{\bot}=A$ and this implies that $u=0$. Hence, $T$ is an $n$-multiplier.
\end{proof}

\end{document}